\documentclass[12pt]{amsart}

\usepackage[all,poly,necula,rotate,tips]{xy}
\usepackage{color}

\usepackage[mathscr]{eucal}
\usepackage{amscd}
\usepackage{amsfonts}
\usepackage{amsmath, amsthm, amssymb}
\usepackage{latexsym}
\usepackage[pdftex]{graphics}

\theoremstyle{plain}
\newtheorem{thm}{Theorem}

\newtheorem{prop}[thm]{Proposition}

\theoremstyle{definition}
\newtheorem{defn}[thm]{Definition}
\newtheorem{ex}[thm]{Example}

\theoremstyle{remark}
\newtheorem{rem}[thm]{Remark}

\numberwithin{equation}{section}

\begin{document}

%% Personnal macros
\def\DeclareMathOperator#1#2{\newcommand#1{\operatorname{#2}}}
\newcommand{\dta}[2]{\Delta^{#1}_{#2}}
\newcommand{\dtb}[2]{{}^{#1}_{#2}}
\newcommand{\ens}[1]
{
 \begin{array}{c}
  \dtb{1}{4}, \dtb{12}{34}, \dtb{123}{234}, \\
  #1
 \end{array}
}

\DeclareMathOperator{\SL}{SL}
\newcommand{\C}{\mathbb{C}}

\makeatletter
\newcommand*{\contains}[2]{%
  \TE@throw
  \noexpand \in@ {#2}{#1}%
  \noexpand \ifin@
}

\newcommand*{\tempb}{}

\newcommand*{\hex}[1]{%
    \renewcommand*{\tempb}{}%
    \@for \tempa :=#1\do {%
    \ifthenelse{\contains{#1}{AC}}%
      {\protected@edef \tempb {\tempb ,"A1";"A3"**\string @{-}}}{}%
    \ifthenelse{\contains{#1}{AD}}%
      {\protected@edef \tempb {\tempb ,"A1";"A4"**\string @{-}}}{}%
    \ifthenelse{\contains{#1}{AE}}
      {\protected@edef \tempb {\tempb ,"A1";"A5"**\string @{-}}}{}%
    \ifthenelse{\contains{#1}{BD}}
      {\protected@edef \tempb {\tempb ,"A2";"A4"**\string @{-}}}{}%
    \ifthenelse{\contains{#1}{BE}}
      {\protected@edef \tempb {\tempb ,"A2";"A5"**\string @{-}}}{}%
    \ifthenelse{\contains{#1}{BF}}
      {\protected@edef \tempb {\tempb ,"A2";"A6"**\string @{-}}}{}%
    \ifthenelse{\contains{#1}{CE}}
      {\protected@edef \tempb {\tempb ,"A3";"A5"**\string @{-}}}{}%
    \ifthenelse{\contains{#1}{CF}}
      {\protected@edef \tempb {\tempb ,"A3";"A6"**\string @{-}}}{}%
    \ifthenelse{\contains{#1}{DF}}
      {\protected@edef \tempb {\tempb ,"A4";"A6"**\string @{-}}}{}%
    }%
    \protected@edef \tempb ##1##2{##1\tempb##2}%
    \tempb {\begin{xy}/l.8pc/:, {\xypolygon6"A"{}}}{\end{xy}}}

\newcommand*{\hexgr}[1]{%
    \renewcommand*{\tempb}{}%
    \@for \tempa :=#1\do {%
    \ifthenelse{\contains{#1}{AC}}%
      {\protected@edef \tempb {\tempb ,"A1";"A3"**\string @{-}}}{}%
    \ifthenelse{\contains{#1}{AD}}%
      {\protected@edef \tempb {\tempb ,"A1";"A4"**\string @{-}}}{}%
    \ifthenelse{\contains{#1}{AE}}
      {\protected@edef \tempb {\tempb ,"A1";"A5"**\string @{-}}}{}%
    \ifthenelse{\contains{#1}{BD}}
      {\protected@edef \tempb {\tempb ,"A2";"A4"**\string @{-}}}{}%
    \ifthenelse{\contains{#1}{BE}}
      {\protected@edef \tempb {\tempb ,"A2";"A5"**\string @{-}}}{}%
    \ifthenelse{\contains{#1}{BF}}
      {\protected@edef \tempb {\tempb ,"A2";"A6"**\string @{-}}}{}%
    \ifthenelse{\contains{#1}{CE}}
      {\protected@edef \tempb {\tempb ,"A3";"A5"**\string @{-}}}{}%
    \ifthenelse{\contains{#1}{CF}}
      {\protected@edef \tempb {\tempb ,"A3";"A6"**\string @{-}}}{}%
    \ifthenelse{\contains{#1}{DF}}
      {\protected@edef \tempb {\tempb ,"A4";"A6"**\string @{-}}}{}%
    }%
    \protected@edef \tempb ##1##2{##1\tempb##2}%
    \tempb {\begin{xy}/l3pc/:, <.3pc, 0pc>, {\xypolygon6"N"{~*{\ifcase\xypolynode\or 1 \or 2 \or 3 \or 4 \or 5 \or 6 \fi }~:{(1.2,0):}~>{}} }, <0pc, 0pc>, {\xypolygon6"A"{}}}{ \end{xy}}}
\makeatother
\newcommand{\vecv}[2]{\left(\!\!
    \begin{array}{c}
      #1 \\
      #2
    \end{array}
    \!\!\right)}
\newcommand{\vech}[2]{\left(\!\!
    \begin{array}{cc}
      #1 & #2
    \end{array}
    \!\!\right)}
\newcommand{\mat}[4]{\left(\!\!
    \begin{array}{cc}
      #1 & #2 \\
      #3 & #4
    \end{array}
    \!\!\right)}
\newcommand{\matt}[9]{\left(\!\!
    \begin{array}{ccc}
      #1 & #2 & #3 \\
      #4 & #5 & #6 \\
      #7 & #8 & #9
    \end{array}
    \!\!\right)}

\newcommand{\dg}{\mathbf{d}}
\newcommand{\N}{\mathbb{N}}
\newcommand{\Z}{\mathbb{Z}}
\DeclareMathOperator{\irr}{irr}
\DeclareMathOperator{\md}{mod}
\DeclareMathOperator{\Gr}{Gr}
\renewcommand{\phi}{\varphi}

\newsavebox{\mybox}
\newlength{\myboxlen}
\newcommand{\figs}[1]
{
 \sbox{\mybox}{#1}
 \settowidth{\myboxlen}{\usebox{\mybox}}
 \begin{minipage}{\myboxlen}
  \usebox{\mybox}
 \end{minipage}
}

\newsavebox{\myboxb}
\newlength{\myboxblen}
\newcommand{\ensb}[5]
{
 \sbox{\myboxb}{\hspace{#1}\fcolorbox{black}{white}{$#3 \oplus #4 \oplus #5$}\hspace{#2}}
 \settowidth{\myboxblen}{\usebox{\myboxb}}
 \begin{minipage}{\myboxblen}
  \usebox{\myboxb}
 \end{minipage}
}

%% reps

\newsavebox{\spu}
\sbox{\spu}
{%
\objectmargin={1pt}$\xymatrix@C=.3cm@R=.3cm{
1 \ar[dr] & & \\
& 2 \ar[dr] & \\
& & 3
}$
}
\newlength{\lpu}
\settowidth{\lpu}{\usebox{\spu}}

\newcommand{\Pu}{\begin{minipage}{\lpu}\usebox{\spu}\end{minipage}}

\newsavebox{\spd}
\sbox{\spd}
{%
\objectmargin={1pt}$\xymatrix@C=.3cm@R=.3cm{
 & 2 \ar[dr] \ar[dl] & \\
1 \ar[dr] & & 3 \ar[dl] \\
& 2 &
}$}
\newlength{\lpd}
\settowidth{\lpd}{\usebox{\spd}}

\newcommand{\Pd}{\begin{minipage}{\lpd}\usebox{\spd}\end{minipage}}

\newsavebox{\spt}
\sbox{\spt}
{%
\objectmargin={1pt}$\xymatrix@C=.3cm@R=.3cm{
 &  & 3 \ar[dl] \\
 & 2 \ar[dl] &  \\
1 & &
}$}
\newlength{\lpt}
\settowidth{\lpt}{\usebox{\spt}}

\newcommand{\Pt}{\begin{minipage}{\lpt}\usebox{\spt}\end{minipage}}

\newsavebox{\ssu}
\sbox{\ssu}
{%
\objectmargin={1pt}$\xymatrix@C=.3cm@R=.3cm{
1
}$
}
\newlength{\lsu}
\settowidth{\lsu}{\usebox{\ssu}}

\newcommand{\Su}{\begin{minipage}{\lsu}\usebox{\ssu}\end{minipage}}

\newsavebox{\ssd}
\sbox{\ssu}
{%
\objectmargin={1pt}$\xymatrix@C=.3cm@R=.3cm{
2
}$
}
\newlength{\lsd}
\settowidth{\lsd}{\usebox{\ssd}}

\newcommand{\Sd}{\begin{minipage}{\lsd}\usebox{\ssd}\end{minipage}}

\newsavebox{\sst}
\sbox{\sst}
{%
\objectmargin={1pt}$\xymatrix@C=.3cm@R=.3cm{
3
}$
}
\newlength{\lst}
\settowidth{\lst}{\usebox{\sst}}

\newcommand{\St}{\begin{minipage}{\lst}\usebox{\sst}\end{minipage}}

\newsavebox{\sdut}
\sbox{\sdut}
{%
\objectmargin={1pt}$\xymatrix@C=.3cm@R=.3cm{
& 2 \ar[dl] \ar[dr] & \\
1 & & 3
}$}
\newlength{\ldut}
\settowidth{\ldut}{\usebox{\sdut}}

\newcommand{\dut}{\begin{minipage}{\ldut}\usebox{\sdut}\end{minipage}}

\newsavebox{\sutd}
\sbox{\sutd}
{%
\objectmargin={1pt}$\xymatrix@C=.3cm@R=.3cm{
1 \ar[dr] & & 3 \ar[dl] \\
& 2 &
}$}
\newlength{\lutd}
\settowidth{\lutd}{\usebox{\sutd}}

\newcommand{\utd}{\begin{minipage}{\lutd}\usebox{\sutd}\end{minipage}}

\newsavebox{\sud}
\sbox{\sud}
{%
\objectmargin={1pt}$\xymatrix@C=.3cm@R=.3cm{
1 \ar[dr] & \\
& 2
}$}
\newlength{\lud}
\settowidth{\lud}{\usebox{\sud}}

\newcommand{\ud}{\begin{minipage}{\lud}\usebox{\sud}\end{minipage}}

\newsavebox{\sdt}
\sbox{\sdt}
{%
\objectmargin={1pt}$\xymatrix@C=.3cm@R=.3cm{
2 \ar[dr] & \\
& 3
}$}
\newlength{\ldt}
\settowidth{\ldt}{\usebox{\sdt}}

\newcommand{\dt}{\begin{minipage}{\ldt}\usebox{\sdt}\end{minipage}}

\newsavebox{\sdu}
\sbox{\sdu}
{%
\objectmargin={1pt}$\xymatrix@C=.3cm@R=.3cm{
& 2 \ar[dl] \\
1 &
}$}
\newlength{\ldu}
\settowidth{\ldu}{\usebox{\sdu}}

\newcommand{\du}{\begin{minipage}{\ldu}\usebox{\sdu}\end{minipage}}

\newsavebox{\std}
\sbox{\std}
{%
\objectmargin={1pt}$\xymatrix@C=.3cm@R=.3cm{
& 3 \ar[dl] \\
2 &
}$}
\newlength{\ltd}
\settowidth{\ltd}{\usebox{\std}}

\newcommand{\td}{\begin{minipage}{\ltd}\usebox{\std}\end{minipage}}

%smallrep

\newsavebox{\spus}
\sbox{\spus}
{%
\tiny\objectmargin={0.3pt}$\xymatrix@C=.2cm@R=.2cm{
1 \ar[dr] & & \\
& 2 \ar[dr] & \\
& & 3
}$
}
\newlength{\lpus}
\settowidth{\lpus}{\usebox{\spus}}

\newcommand{\Pus}{\begin{minipage}{\lpus}\usebox{\spus}\end{minipage}}

\newsavebox{\spds}
\sbox{\spds}
{%
\tiny\objectmargin={0.3pt}$\xymatrix@C=.2cm@R=.2cm{
 & 2 \ar[dr] \ar[dl] & \\
1 \ar[dr] & & 3 \ar[dl] \\
& 2 &
}$}
\newlength{\lpds}
\settowidth{\lpds}{\usebox{\spds}}

\newcommand{\Pds}{\begin{minipage}{\lpds}\usebox{\spds}\end{minipage}}

\newsavebox{\spts}
\sbox{\spts}
{%
\tiny\objectmargin={0.3pt}$\xymatrix@C=.2cm@R=.2cm{
 &  & 3 \ar[dl] \\
 & 2 \ar[dl] &  \\
1 & &
}$}
\newlength{\lpts}
\settowidth{\lpts}{\usebox{\spts}}

\newcommand{\Pts}{\begin{minipage}{\lpts}\usebox{\spts}\end{minipage}}

\newsavebox{\ssus}
\sbox{\ssus}
{%
\tiny\objectmargin={0.3pt}$\xymatrix@C=.2cm@R=.2cm{
1
}$
}
\newlength{\lsus}
\settowidth{\lsus}{\usebox{\ssus}}

\newcommand{\Sus}{\begin{minipage}{\lsus}\usebox{\ssus}\end{minipage}}

\newsavebox{\ssds}
\sbox{\ssus}
{%
\tiny\objectmargin={0.3pt}$\xymatrix@C=.2cm@R=.2cm{
2
}$
}
\newlength{\lsds}
\settowidth{\lsds}{\usebox{\ssds}}

\newcommand{\Sds}{\begin{minipage}{\lsd}\usebox{\ssd}\end{minipage}}

\newsavebox{\ssts}
\sbox{\ssts}
{%
\tiny\objectmargin={0.3pt}$\xymatrix@C=.2cm@R=.2cm{
3
}$
}
\newlength{\lsts}
\settowidth{\lsts}{\usebox{\ssts}}

\newcommand{\Sts}{\begin{minipage}{\lsts}\usebox{\ssts}\end{minipage}}

\newsavebox{\sduts}
\sbox{\sduts}
{%
\tiny\objectmargin={0.3pt}$\xymatrix@C=.2cm@R=.2cm{
& 2 \ar[dl] \ar[dr] & \\
1 & & 3
}$}
\newlength{\lduts}
\settowidth{\lduts}{\usebox{\sduts}}

\newcommand{\duts}{\begin{minipage}{\lduts}\usebox{\sduts}\end{minipage}}

\newsavebox{\sutds}
\sbox{\sutds}
{%
\tiny\objectmargin={0.3pt}$\xymatrix@C=.2cm@R=.2cm{
1 \ar[dr] & & 3 \ar[dl] \\
& 2 &
}$}
\newlength{\lutds}
\settowidth{\lutds}{\usebox{\sutds}}

\newcommand{\utds}{\begin{minipage}{\lutds}\usebox{\sutds}\end{minipage}}

\newsavebox{\suds}
\sbox{\suds}
{%
\tiny\objectmargin={0.3pt}$\xymatrix@C=.2cm@R=.2cm{
1 \ar[dr] & \\
& 2
}$}
\newlength{\luds}
\settowidth{\luds}{\usebox{\suds}}

\newcommand{\uds}{\begin{minipage}{\luds}\usebox{\suds}\end{minipage}}

\newsavebox{\sdts}
\sbox{\sdts}
{%
\tiny\objectmargin={0.3pt}$\xymatrix@C=.2cm@R=.2cm{
2 \ar[dr] & \\
& 3
}$}
\newlength{\ldts}
\settowidth{\ldts}{\usebox{\sdts}}

\newcommand{\dts}{\begin{minipage}{\ldts}\usebox{\sdts}\end{minipage}}

\newsavebox{\sdus}
\sbox{\sdus}
{%
\tiny\objectmargin={0.3pt}$\xymatrix@C=.2cm@R=.2cm{
& 2 \ar[dl] \\
1 &
}$}
\newlength{\ldus}
\settowidth{\ldus}{\usebox{\sdus}}

\newcommand{\dus}{\begin{minipage}{\ldus}\usebox{\sdus}\end{minipage}}

\newsavebox{\stds}
\sbox{\stds}
{%
\tiny\objectmargin={0.3pt}$\xymatrix@C=.2cm@R=.2cm{
& 3 \ar[dl] \\
2 &
}$}
\newlength{\ltds}
\settowidth{\ltds}{\usebox{\stds}}

\newcommand{\tds}{\begin{minipage}{\ltds}\usebox{\stds}\end{minipage}}

\DeclareMathOperator{\Ext}{Ext}
\DeclareMathOperator{\Hom}{Hom}
\DeclareMathOperator{\add}{add}
\DeclareMathOperator{\coker}{coker}

\newcommand{\ensc}[1]
{
 \begin{array}{c}
  \dta{1}{4} = \dta{123}{234}, \dta{12}{34}, \\
  #1
 \end{array}
}

\newcommand{\ag}{\mathbf{a}}

\title{Example of categorification of a cluster algebra}

%% Authors %%%%%%
\author{Laurent Demonet}　% authors　%英文で記入（日本語併記可）

%%%%%%%%%%%%%%　所属先等のデータ　%英文で記入（日本語併記可) %%%%%%%

\address{
\begin{flushleft}
        \hspace{0.3cm}  Graduate School of Mathematics \\
         \hspace{0.3cm}  Nagoya University \\
         \hspace{0.3cm}  Furocho, Chikusaku, Nagoya, 464-8602 Japan \\
\end{flushleft}
}
\email{Laurent.Demonet@normalesup.org}

%%%% 他に投稿予定の方は、必ず下記の内容を脚注に書き加えておいて下さい。%%%%%%%%
%%%%%%%%%%%%%%%%%% 脚注の書式 %%%%%%%%%%%%%%%%%

%% \thanks{The detailed version of this paper will be submitted for publication elsewhere.}

\thanks{The paper is in a final form and no version of it will
 be submitted for publication elsewhere.}

%% \thanks{The detailed /final/ version of this paper will be /has
%% been/ submitted for publication　elsewhere.}

%%%%%%%%%%%%% 本文 %%%%%%%%%%%%%%%%%%%%%

\maketitle
\thispagestyle{empty}

%%%%%%%%%%%%%%%%%% Abstract の書式 （英語）%%%%%%%%%%%%%%%%%

\begin{abstract}

%%% アブストラクト内容 （必ず、英語で記入）

We present here two detailed examples of additive categorifications of the cluster algebra structure of a coordinate ring of a maximal unipotent subgroup of a simple Lie group. The first one is of simply-laced type ($A_3$) and relies on an article by Gei\ss, Leclerc and Schr\"oer. The second is of non simply-laced type ($C_2$) and relies on an article by the author of this note. This is aimed to be accessible, specially for people who are not familiar with this subject.

%%%%%%%%%%% 下記(キーワード、分類)はオプション　%%%%%%%%%%%%%%%
% \bigskip

% {\it Key Words:} \quad Ring, Algebra, Representation.

% \medskip
% {\it $2000$ Mathematics Subject Classification{\rm :}}
% \quad Primary  16Gxx, 16Dxx;  Secondary 16Exx, 16Lxx.
%%%%%%%%%%%%%%%%%%%%%%%%%%%%%%%%%

\end{abstract}

%%%%%%%%%%%%%%%%%% 論文本体 %%%%%%%%%%

\section{Introduction: the total positivity problem}

Let $N$ be the subgroup of $\SL_4(\C)$ consisting of upper triangular matrices with diagonal $1$. We say that $X \in N$ is \emph{totally positive} if its 12 non-trivial minors are positive real numbers (a minor is non-trivial if it is not constant on $N$ and not product of other minors). As a consequence of various results of Fomin and Zelevinsky \cite{FZ1} (see also \cite{BFZ}), in a (very) special case, we get
\begin{prop}[Fomin-Zelevinsky] \label{totpos}
  $X \in N$ is totally positive if and only if the minors $\dta{1}{4}(X)$, $\dta{12}{34}(X)$, $\dta{123}{234}(X)$, $\dta{12}{24}(X)$, $\dta{2}{4}(X)$, $\dta{3}{4}(X)$ are positive.
\end{prop}
where $\dta{\ell_1 \dots \ell_k}{c_1 \dots c_k}(X)$ is the minor of $X$ with rows $\ell_1$, \dots, $\ell_k$ and columns $c_1$, \dots, $c_k$.

Remark that, as the algebraic variety $N$ has dimension $6$, we can not expect to find a criterion with less than $6$ inequalities to check the total positivity of a matrix.

\begin{figure}
 \makeatletter
    \small
    $$\begin{xy}
    		="D",0;<.58mm,.58mm>:
   		{\xypolygon5"O"{~*{\ifcase\xypolynode\or \ens{\dtb{12}{24}, \dtb{2}{4}, \dtb{3}{4}} \or \or \or \or \fi }~:{(0,0):}~>{}}},
   		{\xypolygon5"A"{~*{\ifcase\xypolynode\or \or \ens{\dtb{12}{24}, \dtb{2}{4}, \dtb{12}{23}} \or \ens{\dtb{12}{24}, \dtb{1}{2}, \dtb{3}{4}} \or \or \ens{\dtb{23}{34}, \dtb{2}{4}, \dtb{3}{4}}\fi}~:{(-30,0):}~>{}}},
    		{\xypolygon5"B"{~*{\ifcase\xypolynode\or \ens{\dtb{13}{34}, \dtb{1}{2}, \dtb{3}{4}} \or \ens{\dtb{13}{34}, \dtb{23}{34}, \dtb{3}{4}} \or \ens{\dtb{23}{34}, \dtb{2}{4}, \dtb{2}{3}} \or \ens{\dtb{12}{23}, \dtb{2}{4}, \dtb{2}{3}} \or \ens{\dtb{12}{24}, \dtb{1}{2}, \dtb{12}{23}}\fi}~:{(54,0):}~>{}}},
    		{\xypolygon5"C"{~*{\ifcase\xypolynode\or \ens{\dtb{13}{34}, \dtb{1}{2}, \dtb{1}{3}} \or \ens{\dtb{13}{34}, \dtb{23}{34}, \dtb{1}{3}} \or \ens{\dtb{23}{34}, \dtb{2}{3}, \dtb{1}{3}} \or \ens{\dtb{12}{23}, \dtb{2}{3}, \dtb{1}{3}} \or \ens{\dtb{12}{23}, \dtb{1}{2}, \dtb{1}{3}}\fi}~:{(90,0):}~>{}}},
    		"B1";"C1"**\dir{-},
    		"B2";"C2"**\dir{-},
    		"B3";"C3"**\dir{-},
    		"B4";"C4"**\dir{-},
    		"B5";"C5"**\dir{-},
    		"B1";"A3"**\dir{-},
    		"A3";"B5"**\dir{-},
    		"B5";"A2"**\dir{-},
    		"A2";"B4"**\dir{-},
    		"B4";"B3"**\dir{-},
    		"B3";"A5"**\dir{-},
    		"A5";"B2"**\dir{-},
    		"B2";"B1"**\dir{-},
    		"A3";"O1"**\dir{-},
    		"A2";"O1"**\dir{-},
    		"A5";"O1"**\dir{-},
    		"C1";"C2"**\dir{-},
    		"C2";"C3"**\dir{-},
    		"C3";"C4"**\dir{-},
    		"C4";"C5"**\dir{-},
    		"C5";"C1"**\dir{-}
    \end{xy}$$
\makeatother
 \caption{Exchange graph of minors} \label{assoc}
\end{figure}
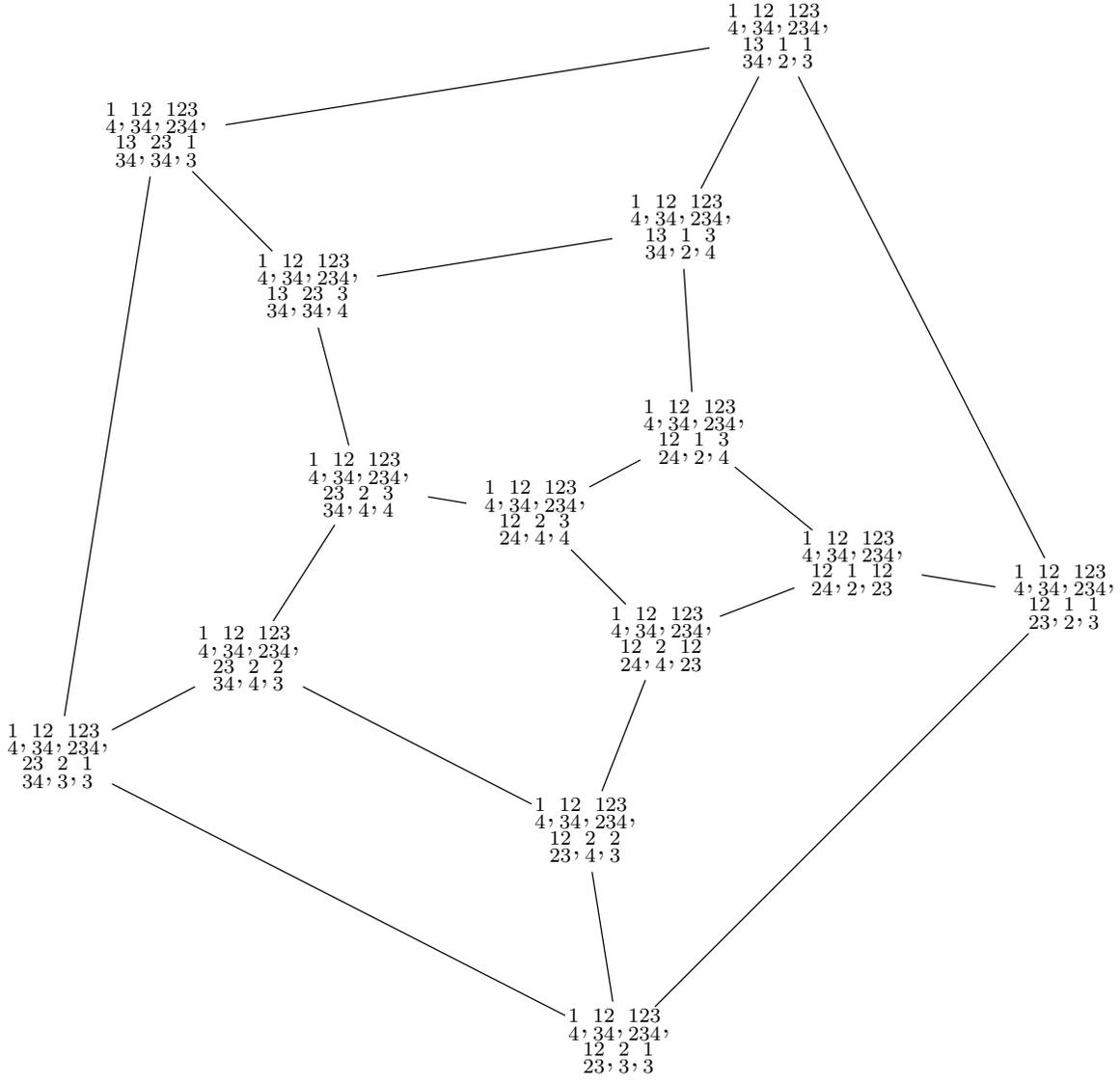

To prove this, just remark that we have the following equality:
$$\dta{12}{24} \dta{23}{34} = \dta{123}{234} \dta{2}{4} + \dta{3}{4} \dta{12}{34}$$
which immediately implies that $\dta{1}{4}(X)$, $\dta{12}{34}(X)$, $\dta{123}{234}(X)$, $\dta{12}{24}(X)$, $\dta{2}{4}(X)$, $\dta{3}{4}(X)$ are positive if and only if $\dta{1}{4}(X)$, $\dta{12}{34}(X)$, $\dta{123}{234}(X)$, $\dta{23}{34}(X)$, $\dta{2}{4}(X)$, $\dta{3}{4}(X)$ are positive. Such an equality is called an \emph{exchange identity}. In Figure \ref{assoc}, we wrote $14$ sets of minors which are related by exchange identities whenever they are linked by an edge. As every minor appears in this graph, it induces the previous proposition.

These observations lead to the definition of a \emph{cluster algebra} \cite{FZ2}. A cluster algebra is an algebra endowed with an additional combinatorial structure. Namely, a (generally infinite) set of distinguished elements called \emph{cluster variables} grouped into subsets of the same cardinality $n$, called \emph{clusters} and a finite set $\{x_{n+1}, x_{n+2}, \dots, x_m\}$ called the set of \emph{coefficients}. For each cluster $\{x_1, x_2, \dots, x_n\}$, the \emph{extended cluster} $\{x_1, \dots, x_n, x_{n+1}, \dots, x_m\}$ is a transcendence basis of the algebra. Moreover, each cluster $\{x_1, x_2, \dots, x_n\}$ has $n$ neighbours obtained by replacing one of its elements $x_k$ by a new one $x'_k$ related by a relation
$$x_k x'_k = M_1 + M_2$$
where $M_1$ and $M_2$ are mutually prime monomials in $\{x_1, \dots, x_{k-1}, x_{k+1}, \dots, x_m\}$, given by precise combinatorial rules. These replacements, called \emph{mutations} and denoted by $\mu_k$ are involutive. For precise definitions and details about these constructions, we refer to \cite{FZ2}.

In the previous example, the coefficients are $\dta{1}{4}$, $\dta{12}{34}$ and $\dta{123}{234}$ and the cluster variables are all the other non-trivial minors. The extended clusters are the sets appearing at the vertices of Figure \ref{assoc}.

The aim of the following sections is to describe examples of \emph{additive categorifications} of cluster algebras. It consists of enhancing the cluster algebra structure with an additive category, some objects of which reflect the combinatorial structure of the cluster algebra; moreover, there is an explicit formula, the \emph{cluster character} associating to these particular objects elements of the algebra, in a way which is compatible with the combinatorial structure. The examples we develop here rely on (abelian) module categories. They are particular cases of categorifications by exact categories appearing in \cite{GLS} (simply-laced case) and \cite{Dem} (non simply-laced case). The study of cluster algebras and their categorifications has been particularly successful these last years. For a survey on categorification by triangulated categories and a much more complete bibliography, see \cite{Kel1}.

\section{The preprojective algebra and the cluster character}

Let $Q$ be the following quiver (oriented graph):
$$\xymatrix{
 1 \ar@/^/[r]^\alpha & 2 \ar@/_/[r]_{\beta^*} \ar@/^/[l]^{\alpha^*} & 3 \ar@/_/[l]_{\beta}
}$$

As usual, denote by $\C Q$ the $\C$-algebra, a basis of which is formed by the paths (including $0$-length paths supported by each of the three vertices) and the multiplication of which is defined by concatenation of paths when it is possible and vanishes when paths can not be composed (we write here the composition from left to right, on the contrary to the usual composition of maps). Thus, a (right) $\C Q$-module is naturally graded by idempotents ($0$-length paths) corresponding to vertices and the action of arrows seen as elements of the algebra can naturally be identified with linear maps between the corresponding homogeneous subspaces of the representation. We shall use the following right-hand side convenient notation:
 $$\xymatrix{
  \C \ar@/^/[r]^{\vecv{0}{-1}} & \C^2 \ar@/_/[r]_{\mat{1}{0}{0}{1}} \ar@/^/[l]^{\vech{1}{0}} & \C^2 \ar@/_/[l]_{\mat{0}{0}{1}{0}}
 } = 
 \figs{\objectmargin={1pt}\xymatrix@C=.3cm@R=.3cm{
  & 2 \ar[dl] \ar[dr] & \\
  1 \ar[dr]_{-1} & & 3 \ar[dl] \\
  & 2 \ar[dr] & \\
  & & 3
 }}$$
where each of the digits represents a basis vector of the representation and each arrow a non-zero scalar ($1$ when not specified) in the corresponding matrix entry.

Let us now introduce the preprojective algebra of $Q$:

\begin{defn}
 The \emph{preprojective algebra} of $Q$ is defined by
 $$\Pi_Q = \frac{\C Q}{\left(\alpha \alpha^*, \alpha^* \alpha + \beta^* \beta, \beta \beta^* \right)}$$
\end{defn}
the representations of which are seen as particular representations of $\C Q$ (in other words, $\md \Pi_Q$ is a full subcategory of $\md \C Q$).

\begin{ex}
 Among the following representations of $\C Q$, the first one and the second one are representations of $\Pi_Q$:
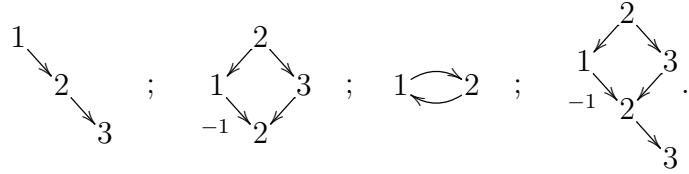
 $$
 \figs{\objectmargin={1pt}\xymatrix@C=.3cm@R=.3cm{
  1 \ar[dr] & & \\
  & 2 \ar[dr] & \\
  & & 3
 }} \quad ; \quad
 \figs{\objectmargin={1pt}\xymatrix@C=.3cm@R=.3cm{
  & 2 \ar[dl] \ar[dr] & \\
  1 \ar[dr]_{-1} & & 3 \ar[dl] \\
  & 2 & \\
 }} \quad ; \quad
 \figs{\objectmargin={1pt}\xymatrix@C=.3cm@R=.3cm{
  1 \ar@/^/[rr] & & 2 \ar@/^/[ll]
 }} \quad ; \quad
 \figs{\objectmargin={1pt}\xymatrix@C=.3cm@R=.3cm{
  & 2 \ar[dl] \ar[dr] & \\
  1 \ar[dr]_{-1} & & 3 \ar[dl] \\
  & 2 \ar[dr] & \\
  & & 3
 }}.$$
\end{ex}

One of the property, which is discussed in many places (for example in \cite{GLS}), of the preprojective algebra of $Q$, fundamental for this categorification, is

\begin{prop}
 The category $\md \Pi_Q$ is \emph{stably $2$-Calabi-Yau}. In other words, for every $X, Y \in \md \Pi_Q$,
 $$\Ext^1(X, Y) \simeq \Ext^1(Y,X)^*$$
 functorially in $X$ and $Y$, where $\Ext^1(Y, X)^*$ is the $\C$-dual of $\Ext^1(Y,X)$. In particular, it is a Frobenius category (is has enough projective objects and enough injective objects and they coincide).
\end{prop}

Let us now define the three following one-parameter subgroups of $N$:
$$x_1(t) = \left( \begin{matrix} 1 & t & 0 & 0 \\ 0 & 1 & 0 & 0 \\  0 & 0 & 1 & 0 \\ 0 & 0 & 0 & 1 \end{matrix} \right) \quad x_2(t) = \left( \begin{matrix} 1 & 0 & 0 & 0 \\ 0 & 1 & t & 0 \\  0 & 0 & 1 & 0 \\ 0 & 0 & 0 & 1 \end{matrix} \right) \quad  x_3(t) = \left( \begin{matrix} 1 & 0 & 0 & 0 \\ 0 & 1 & 0 & 0 \\  0 & 0 & 1 & t \\ 0 & 0 & 0 & 1 \end{matrix} \right).$$

For $X \in \md \Pi_Q$ and any sequence of vertices $a_1, a_2, \dots, a_n$ of $Q$, we denote by
$$\Phi_{X, a_1 a_2 \dots a_n} = \left\{0 = X_0 \subset X_1 \subset \dots \subset X_{n-1} \subset X_n = X \,|\, \forall i \in \{1, 2, \dots, n\}, \frac{X_i}{X_{i-1}} \simeq S_{a_i}\right\}$$
the \emph{variety of composition series} of $X$ of type $a_1 a_2 \dots a_n$ ($S_{a_i}$ is the simple module, of dimension $1$, supported at vertex $a_i$). This is a closed algebraic subvariety of the product of Grassmannians
$$\Gr_1(X) \times \Gr_2(X) \times \dots \times \Gr_n(X).$$
We denote by $\chi$ the Euler characteristic. Using results of Lusztig and Kashiwara-Saito, Gei\ss-Leclerc-Schro\"er proved the following result:

\begin{thm}[\cite{GLS}]
  Let $X \in \md \Pi_Q$. There is a unique $\phi_X \in \C[N]$ such that

  $$ \phi_X\left(x_{a_1}(t_1)x_{a_2}(t_2) \dots x_{a_6}(t_6)\right) = \sum_{i_1, i_2, \dots, i_6 \in \N} \chi\left(\Phi_{X, a_1^{i_1} a_2^{i_2} \dots a_6^{i_6}} \right) \frac{t_1^{i_1} t_2^{i_2} \dots t_6^{i_6}}{i_1!i_2!\dots i_6!}$$

  for every word $a_1 a_2 a_3 a_4 a_5 a_6$ representing the longest element of $\mathfrak{S}_4$ ($a_k^{i_k}$ is the repetition $i_k$ times of $a_k$).
 \end{thm}
The map $\phi : \md \Pi_Q \rightarrow \C[N]$ is called a \emph{cluster character}.

\begin{rem}
 \begin{enumerate}
  \item The uniqueness in the previous theorem is easy because it is well known that
   $$x_{a_1}(t_1)x_{a_2}(t_2) \dots x_{a_6}(t_6)$$
   runs over a dense subset of $N$ ;
  \item the existence is much harder and strongly relies on the construction of semi-canonical bases by Lusztig \cite{L}. In particular, the fact that it does not depend on the choice of $a_1 a_2 a_3 a_4 a_5 a_6$ is not clear \emph{a priori} (see the following examples).
 \end{enumerate}
\end{rem}

\begin{ex}
 We suppose that $a_1 a_2 a_3 a_4 a_5 a_6 = 213213$. Then
 $$x_{a_1}(t_1)x_{a_2}(t_2) x_{a_3}(t_3)x_{a_4}(t_4)x_{a_5}(t_5) x_{a_6}(t_6) = \left( \begin{matrix} 1 & t_2+t_5 & t_2 t_4 & t_2 t_4 t_6 \\ 0 & 1 & t_1+t_4 & t_1 t_3 + t_1 t_6 + t_4 t_6 \\  0 & 0 & 1 & t_3+t_6 \\ 0 & 0 & 0 & 1 \end{matrix} \right).$$
 \begin{itemize}
  \item The module $S_1$ has only one composition series, of type $1$. Therefore $\Phi_1(S_1)$ is one point and $\Phi_{\ag}(S_1) = \emptyset$ for any other $\ag$. Identifying the two members in the formula of the previous theorem,
   $$\phi_{S_1}\left(x_{a_1}(t_1)x_{a_2}(t_2) x_{a_3}(t_3)x_{a_4}(t_4)x_{a_5}(t_5) x_{a_6}(t_6)\right) = t_2 + t_5 = \dta{1}{2}.$$
  \item The module
   $$P_2 = \Pd$$
   has two composition series, of type $2312$ and $2132$. Therefore,
   $$\phi_{P_2}\left(x_{a_1}(t_1)x_{a_2}(t_2) x_{a_3}(t_3)x_{a_4}(t_4)x_{a_5}(t_5) x_{a_6}(t_6)\right) = t_1 t_2 t_3 t_4 = \dta{12}{34}.$$
   Remark that, in this case, the only composition series which is playing a role is $2132$, even if the situation is symmetric. This justify the second part of the previous remark.
 \end{itemize}
 The other indecomposable representations of $\Pi_Q$ and their cluster character values are collected in Figure \ref{char}.
\end{ex}

\begin{figure}
 \begin{center}
  \begin{tabular}{|c|c|c|c|c|c|c|c|}
   \hline
    \rule[-3ex]{0pt}{7ex} $X \in \md \Pi_Q$ &
    $S_1$ &
    $S_2$ &
    $S_3$ &
    $\ud$ &
    $\du$ &
    $\dt$ &
    $\td$ \\
   \hline
    \rule[-1ex]{0pt}{3.5ex} $\phi_X \in \C[N]$
    & $\dta{1}{2}$
    & $\dta{2}{3}$
    & $\dta{3}{4}$
    & $\dta{12}{23}$
    & $\dta{1}{3}$
    & $\dta{23}{34}$
    & $\dta{2}{4}$ \\
   \hline
  \end{tabular}
  \begin{tabular}{|c|c|c|c|c|c|}
   \hline
    \rule[-4.5ex]{0pt}{10ex} $X \in \md \Pi_Q$ &
    $\dut$ &
    $\utd$ &
    $\Pu$ &
    $\Pd$ &
    $\Pt$ \\
   \hline
    \rule[-1ex]{0pt}{3.5ex} $\phi_X \in \C[N]$
    & $\dta{13}{34}$
    & $\dta{12}{24}$
    & $\dta{123}{234}$
    & $\dta{12}{34}$
    & $\dta{1}{4}$ \\
   \hline
  \end{tabular}
 \end{center}
 \caption{Cluster character} \label{char}
\end{figure}

Two important properties of this cluster character were proved by Gei\ss-Leclerc-Schro\"er (see for example \cite{GLS}):
\begin{prop} \label{charprop}
 Let $X, Y \in \md \Pi_Q$.
 \begin{enumerate}
  \item $\phi_{X \oplus Y} = \phi_X \phi_Y$.
  \item Suppose that $\dim \Ext^1(X, Y) = 1$ (and therefore $\dim \Ext^1(Y, X) = 1$) and let
  $$0 \rightarrow X \rightarrow T_a \rightarrow Y \rightarrow 0 \quad \text{and} \quad 0 \rightarrow Y \rightarrow T_b \rightarrow X \rightarrow 0$$
  be two (unique up to isomorphism) non-split short exact sequences. Then
  $$\phi_X \phi_Y = \phi_{T_a} + \phi_{T_b}.$$
 \end{enumerate}
\end{prop}

\section{Minimal approximations}

This section recall the definition and elementary properties of approximations. It is there for the sake of ease. In what follows, $\md \Pi_Q$ can be replaced by any additive $\Hom$-finite category over a field.

\begin{defn}
 Let $X$ and $T$ be two objects of $\md \Pi_Q$. A \emph{left $\add(T)$-approximation} of $X$ is a morphism $f:X \rightarrow T'$ such that
 \begin{itemize}
  \item $T' \in \add(T)$ (which means that every indecomposable summand of $T'$ is an indecomposable summand of $T$) ;
  \item every morphism $g: X \rightarrow T$ factors through $f$.
 \end{itemize}

 If, moreover, there is no strict direct summand $T''$ of $T'$ and left $add(T)$-approximation $f':X \rightarrow T''$, then $f$ is said to be a \emph{minimal left $\add(T)$-approximation}.
\end{defn}

In the same way, we can define

\begin{defn}
 Let $X$ and $T$ be two objects in $\md \Pi_Q$. A \emph{right $\add(T)$-approximation} of $X$ is a morphism $f:T' \rightarrow X$ such that
 \begin{itemize}
  \item $T' \in \add(T)$ ;
  \item every morphism $g: T \rightarrow X$ factors through $f$.
 \end{itemize}

 If, moreover, there is no strict direct summand $T''$ of $T'$ and right $add(T)$-approximation $f':T'' \rightarrow X$, then $f$ is said to be a \emph{minimal right $\add(T)$-approximation}.
\end{defn}

Now, a classical proposition which permits to explicitly compute approximations:

\begin{prop} \label{compapp}
 Let $X$ and $T \simeq T_1^{i_1} \oplus T_2^{i_2} \oplus \dots \oplus T_n^{i_n}$ be two objects in $\md \Pi_Q$ (the $T_i$'s are non-isomorphic indecomposable). For $i, j \in \{1, \dots, n\}$, we denote by $I_{ij}$ the subvector space of $\Hom(T_i, T_j)$ consisting of the non-invertible morphisms ($I_{ij} = \Hom(T_i, T_j)$ if $i \neq j$). Thus, for $j \in \{1, \dots, n\}$, we obtain a linear map
 \begin{align*}
  \bigoplus_{i \in \{1, \dots, n\}} I_{ij} \otimes \Hom(X, T_i) & \xrightarrow{\phi_j} \Hom(X, T_j) \\
  (g, f) & \mapsto g \circ f.
 \end{align*}
 Let $\mathcal{B}_j$ be a basis of $\coker \phi_j$ lifted to $\Hom(X, T_j)$. Then the morphism
 $$X \xrightarrow{\left(f\right)_{j \in \{1, \dots, n\}, f \in \mathcal{B}_j}} \bigoplus_{j \in \{1, \dots, n\}} T_j^{\# \mathcal{B}_j}$$
 is a minimal left $\add(T)$-approximation of $X$. Moreover, any minimal left $\add(T)$-approximation of $X$ is isomorphic to it.
\end{prop}

The previous proposition has a dual version which permits to compute minimal right approximations. In practice, this computation relies on searching morphisms up to factorization through other objects. There is an explicit example of computation in Example \ref{exmut}.

\section{Maximal rigid objects and their mutations}

Let us introduce the objects the combinatorics of which will play the role of the cluster algebra structure.

\begin{defn}
 Let $X \in \md \Pi_Q$.
 \begin{itemize}
  \item The module $X$ is said to be \emph{rigid} if it has no self-extension, (\emph{i.e.}, $\Ext^1(X, X) = 0$).
  \item The module $X$ is said to be \emph{basic maximal rigid} if it is basic (\emph{i.e.}, it does not have two isomorphic indecomposable summands), rigid, and maximal for these two properties.
 \end{itemize}
\end{defn}

\begin{rem}
 A basic maximal rigid $\Pi_Q$-module contains $\Pi_Q$ as a direct summand (because $\Pi_Q$ is both projective and injective and therefore has no extension with any module).
\end{rem}

\begin{ex}
 The object
 $$\utd \oplus \td \oplus \ud \oplus \Pu \oplus \Pd \oplus \Pt,$$
 the last three summands of which are the indecomposable projective-injective $\Pi_Q$-modules, is basic maximal rigid. It is easy to check that it is basic and rigid, but more difficult to prove that it is maximal for these properties (see \cite{GLS} for more details).
\end{ex}

\begin{rem}
 We can prove that all basic maximal rigid objects have the same number of indecomposable summands (six in the example we are talking about).
\end{rem}

The following result permits to define a mutation on basic maximal rigid objects. Considered as an operation on isomorphism classes of basic maximal rigid objects, the induced combinatorial structure will correspond to the one of a cluster algebra.

\begin{thm}[\cite{GLS}] \label{thmut}
  Let $T \simeq T_1 \oplus T_2 \oplus T_3 \oplus P_1 \oplus P_2 \oplus P_3  \in \md \Pi_Q$ be basic maximal rigid such that $P_1$, $P_2$ and $P_3$ are the indecomposable projective $\Pi_Q$-modules and $T_1$, $T_2$ and $T_3$ are indecomposable non-projective $\Pi_Q$-modules. Then, for $i \in \{1, 2, 3\}$, there are two (unique) short exact sequences
  $$0 \rightarrow T_i \xrightarrow{f} T_a \xrightarrow{f'} T_i^* \rightarrow 0 \quad \text{and} \quad
  0 \rightarrow T_i^* \xrightarrow{g} T_b \xrightarrow{g'} T_i \rightarrow 0$$
  such that
  \begin{enumerate}
   \item $f$ and $g$ are minimal left $\add(T/T_i)$-approximations ;
   \item $f'$ and $g'$ are minimal right $\add(T/T_i)$-approximations ;
   \item $T_i^*$ is indecomposable and non-projective ;
   \item $\dim \Ext^1(T_i, T_i^*) = \dim \Ext^1(T_i^*, T_i) = 1$ and the two short exact sequences do not split ;
   \item $\mu_i(T) = T / T_i \oplus T_i^*$ is basic maximal rigid ;
   \item $T_a$ and $T_b$ do not have common summands.
  \end{enumerate}
\end{thm}

\begin{rem}
 In the previous theorem, the existence and uniqueness, regarding the first two conditions, are automatic, except the fact that the extremities of the two short exact sequences coincide up to order. This fact strongly relies on the stably $2$-Calabi-Yau property. It implies that $\mu_i$ is involutive.
\end{rem}

\begin{defn}
 In the previous theorem, $\mu_i$ is called the \emph{mutation in direction $i$}. The short exact sequences appearing are called \emph{exchange sequences}.
\end{defn}

\begin{ex}
 \label{exmut}
 Let
 $$T = \utd \oplus \td \oplus \ud \oplus \Pu \oplus \Pd \oplus \Pt.$$
 Using Proposition \ref{compapp}, we get a left $\add\left(T / \td\right)$-approximation of $\td$:
 $$\td \rightarrow \utd$$
 and computing the cokernel, we get the exchange sequence:
 $$0 \rightarrow \td \rightarrow \utd \rightarrow S_1 \rightarrow 0$$
 so that
 $$\mu_2(T) = \utd \oplus S_1 \oplus \ud \oplus \Pu \oplus \Pd \oplus \Pt.$$
 Doing mutation in the reverse direction:
 $$0 \rightarrow S_1 \rightarrow \Pt \rightarrow \td \rightarrow 0.$$
 Let us now compute $\mu_1\mu_2(T)$ with its two exchange sequences:
 $$0 \rightarrow \utd \rightarrow S_1 \oplus \Pd \rightarrow \du \rightarrow 0$$
 $$0 \rightarrow \du \rightarrow \ud \oplus \Pt \rightarrow \utd \rightarrow 0$$
 $$\mu_1\mu_2(T) = \du \oplus S_1 \oplus \ud \oplus \Pu \oplus \Pd \oplus \Pt.$$
\end{ex}

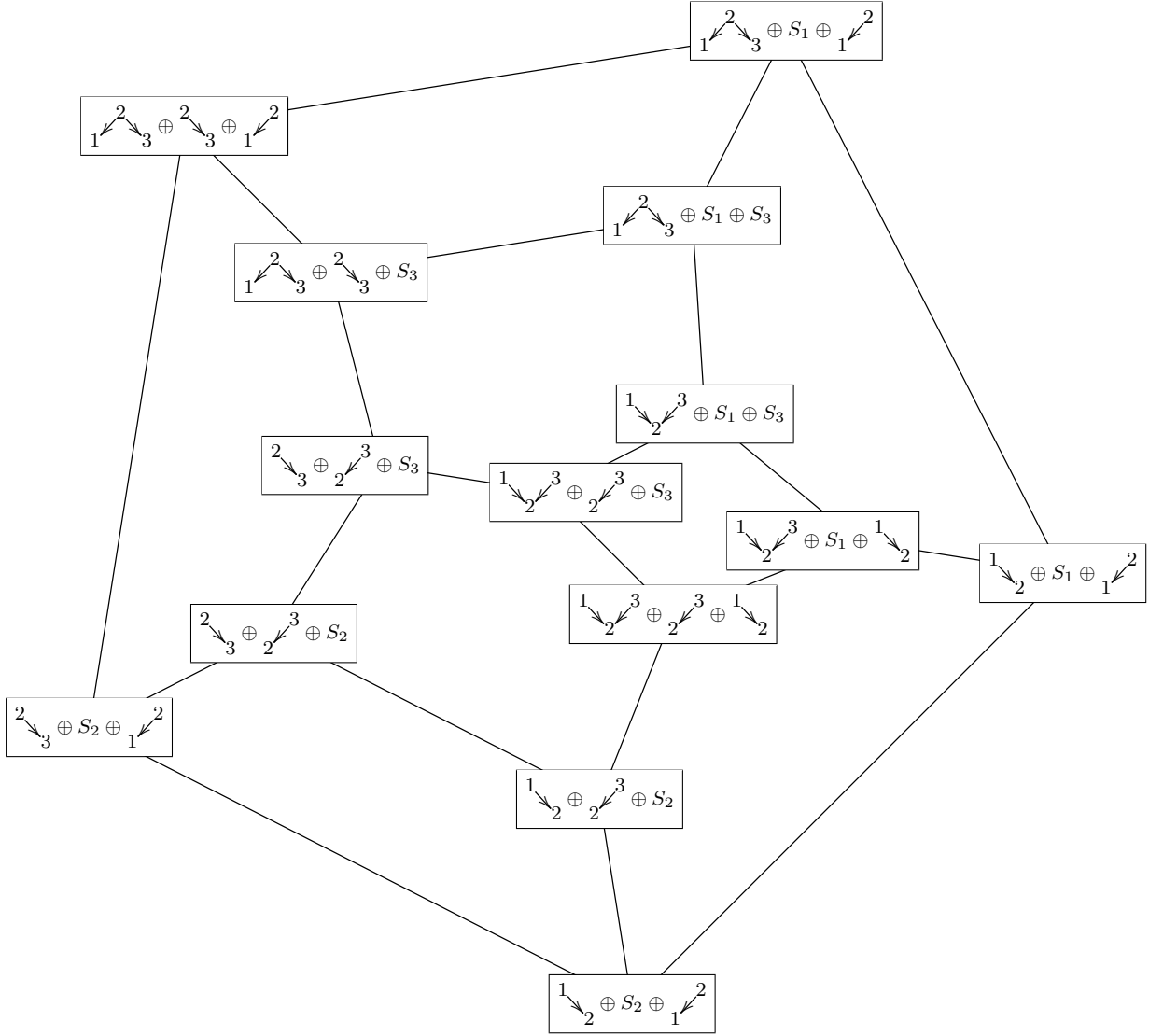
\begin{figure}
  \makeatletter
    \tiny  \objectmargin={-1pt}
    $$\begin{xy}
    		="D",0;<.58mm,.58mm>:
   		{\xypolygon5"O"{~:{(0,0):}~>{}}},
   		{\xypolygon5"A"{~:{(-30,0):}~>{}}},
    		{\xypolygon5"B"{~:{(54,0):}~>{}}},
    		{\xypolygon5"C"{~:{(90,0):}~>{}}},
    		"B1";"C1"**\dir{-},
    		"B2";"C2"**\dir{-},
    		"B3";"C3"**\dir{-},
    		"B4";"C4"**\dir{-},
    		"B5";"C5"**\dir{-},
    		"B1";"A3"**\dir{-},
    		"A3";"B5"**\dir{-},
    		"B5";"A2"**\dir{-},
    		"A2";"B4"**\dir{-},
    		"B4";"B3"**\dir{-},
    		"B3";"A5"**\dir{-},
    		"A5";"B2"**\dir{-},
    		"B2";"B1"**\dir{-},
    		"A3";"O1"**\dir{-},
    		"A2";"O1"**\dir{-},
    		"A5";"O1"**\dir{-},
    		"C1";"C2"**\dir{-},
    		"C2";"C3"**\dir{-},
    		"C3";"C4"**\dir{-},
    		"C4";"C5"**\dir{-},
    		"C5";"C1"**\dir{-},
                "O1"*{\ensb{10mm}{0mm}{\utds}{\tds}{S_3}},
		"A2"*{\ensb{0mm}{0mm}{\utds}{\tds}{\uds}},
		"A3"*{\ensb{0mm}{0mm}{\utds}{S_1}{S_3}},
		"A5"*{\ensb{0mm}{10mm}{\dts}{\tds}{S_3}},
		"B1"*{\ensb{0mm}{0mm}{\duts}{S_1}{S_3}},
		"B2"*{\ensb{0mm}{0mm}{\duts}{\dts}{S_3}},
		"B3"*{\ensb{0mm}{0mm}{\dts}{\tds}{S_2}},
		"B4"*{\ensb{0mm}{0mm}{\uds}{\tds}{S_2}},
		"B5"*{\ensb{0mm}{10mm}{\utds}{S_1}{\uds}},
		"C1"*{\ensb{0mm}{0mm}{\duts}{S_1}{\dus}},
		"C2"*{\ensb{0mm}{0mm}{\duts}{\dts}{\dus}},
		"C3"*{\ensb{0mm}{0mm}{\dts}{S_2}{\dus}},
		"C4"*{\ensb{0mm}{0mm}{\uds}{S_2}{\dus}},
		"C5"*{\ensb{0mm}{0mm}{\uds}{S_1}{\dus}},
    \end{xy}$$
\makeatother
 \caption{Exchange graph of maximal rigid objects (up to projective summands)} \label{assocrep}
\end{figure}

Computing inductively all the mutations, we obtain the \emph{exchange graph of maximal rigid objects of $\Pi_Q$} (Figure \ref{assocrep}).

Then, using Proposition \ref{charprop} and Theorem \ref{thmut} together with other technical results, we get the following proposition:

\begin{prop}[\cite{GLS}]
 If we project the mutation of maximal rigid objects to $\C[N]$ through the cluster character $\phi$, we get a cluster algebra structure on $\C[N]$ (in the sense of \cite{FZ2}). Moreover, this structure is the one proposed combinatorially in \cite{BFZ}. Under this projection, we get the correspondence:
 \begin{align*}
  \{\text{non projective indecomposable objects}\} &\leftrightarrow \{\text{cluster variables}\} \\
  \{\text{projective indecomposable objects}\} &\leftrightarrow \{\text{coefficients}\} \\
  \{\text{basic maximal rigid objects}\} &\leftrightarrow \{\text{extended clusters}\}
 \end{align*}

\end{prop}

\begin{ex}
 Taking the notation of Example \ref{exmut} and looking at Figure \ref{char}, we get:
 $$\dta{1}{2} \dta{2}{4} = \phi_{S_1} \phi_{\td} = \phi_{\utd} + \phi_{\Pt} = \dta{12}{24} + \dta{1}{4}$$
 and
 \begin{align*} \dta{12}{24} \dta{1}{3} &= \phi_{\utd} \phi_{\du} = \phi_{S_1 \oplus \Pd} + \phi_{\ud \oplus \Pt} \\ &= \phi_{S_1} \phi_{\Pd} + \phi_{\ud} \phi_{\Pt} = \dta{1}{2} \dta{12}{34} + \dta{12}{23} \dta{1}{4}.\end{align*}
 which can be easily checked by hand. These are part of the equalities which appear in the proof of Proposition \ref{totpos}.
\end{ex}

\section{From simply-laced case to general one}

Define the following symplectic form:
$$\Psi = \left( \begin{matrix} 0 & 0 & 0 & 1 \\ 0 & 0 & -1 & 0 \\ 0 & 1 & 0 & 0 \\ -1 & 0 & 0 & 0 \end{matrix} \right)$$
and the subgroup
$$N' = \{M \in N| {}^t M \Psi M = \Psi\} \quad \text{or, equivalently} \quad N' = N^{\Z/2\Z}$$
where $\Z/2\Z = \langle g \rangle$ acts on $N$ by $M \mapsto \Psi^{-1}\left({}^t M^{-1}\right) \Psi$. The group $N'$ is a maximal unipotent subgroup of a symplectic group of type $C_2$.

The only non-trivial action of $\Z/2\Z$ on $Q$ induces an action on $\Pi_Q$ and therefore on $\md \Pi_Q$. Denote by $\pi : \C[N] \rightarrow \C[N']$ the canonical projection. We can now formulate the following result:

\begin{thm}[\cite{Dem}]
 \begin{enumerate}
  \item If $T$ is a $\Z/2\Z$-stable basic maximal rigid $\Pi_Q$-module, then $\mu_1\mu_3(T) = \mu_3\mu_1(T)$. Moreover, $\mu_1\mu_3(T)$ and $\mu_2(T)$ are also $\Z/2\Z$-stable.
  \item If $X \in \md \Pi_Q$, then $\pi\left(\phi_X\right) = \pi\left(\phi_{g X}\right)$.
  \item If we denote $\bar \mu_2 = \mu_2$ and $\bar \mu_1 = \mu_1 \mu_3 = \mu_3 \mu_1$, acting on the set of $\Z/2\Z$-stable maximal rigid $\Pi_Q$-modules, $\bar \mu$ induces through $\pi \circ \phi$ the structure of a cluster algebra on $\C[N']$, the clusters of which are projections of the $\Z/2\Z$-stable ones of $\C[N]$.
 \end{enumerate}
\end{thm}

\begin{ex}
 We have
 $$\dta{12}{23}\left( \begin{matrix} 1 & a_{12} & a_{13} & a_{14} \\ 0 & 1 & a_{23} & a_{24} \\ 0 & 0 & 1 & a_{34} \\ 0 & 0 & 0 & 1 \end{matrix} \right) = a_{12} a_{23} - a_{13} \quad \text{and} \quad \dta{2}{4}\left( \begin{matrix} 1 & a_{12} & a_{13} & a_{14} \\ 0 & 1 & a_{23} & a_{24} \\ 0 & 0 & 1 & a_{34} \\ 0 & 0 & 0 & 1 \end{matrix} \right) = a_{24}.$$
 Moreover,
 \begin{align*} & \Psi^{-1}{\vphantom{\begin{matrix} 1 & a_{12} & a_{13} & a_{14} \\ 0 & 1 & a_{23} & a_{24} \\ 0 & 0 & 1 & a_{34} \\ 0 & 0 & 0 & 1 \end{matrix}}}^t\left(\begin{matrix} 1 & a_{12} & a_{13} & a_{14} \\ 0 & 1 & a_{23} & a_{24} \\ 0 & 0 & 1 & a_{34} \\ 0 & 0 & 0 & 1 \end{matrix} \right)^{-1} \Psi \\= & \left(\begin{matrix} 1 & a_{34} & a_{23} a_{34} - a_{24}& a_{12} a_{23} a_{34} - a_{12} a_{24} - a_{13} a_{34} + a_{14} \\ 0 & 1 & a_{23} & a_{12} a_{23} - a_{13} \\ 0 & 0 & 1 & a_{12} \\ 0 & 0 & 0  & 1 \end{matrix} \right) \end{align*}
 which implies that, as expected, $$\pi\left(\dta{12}{23}\right) = \pi\phi_{\ud} = \pi\phi_{\td} = \pi\left(\dta{2}{4}\right).$$
\end{ex}

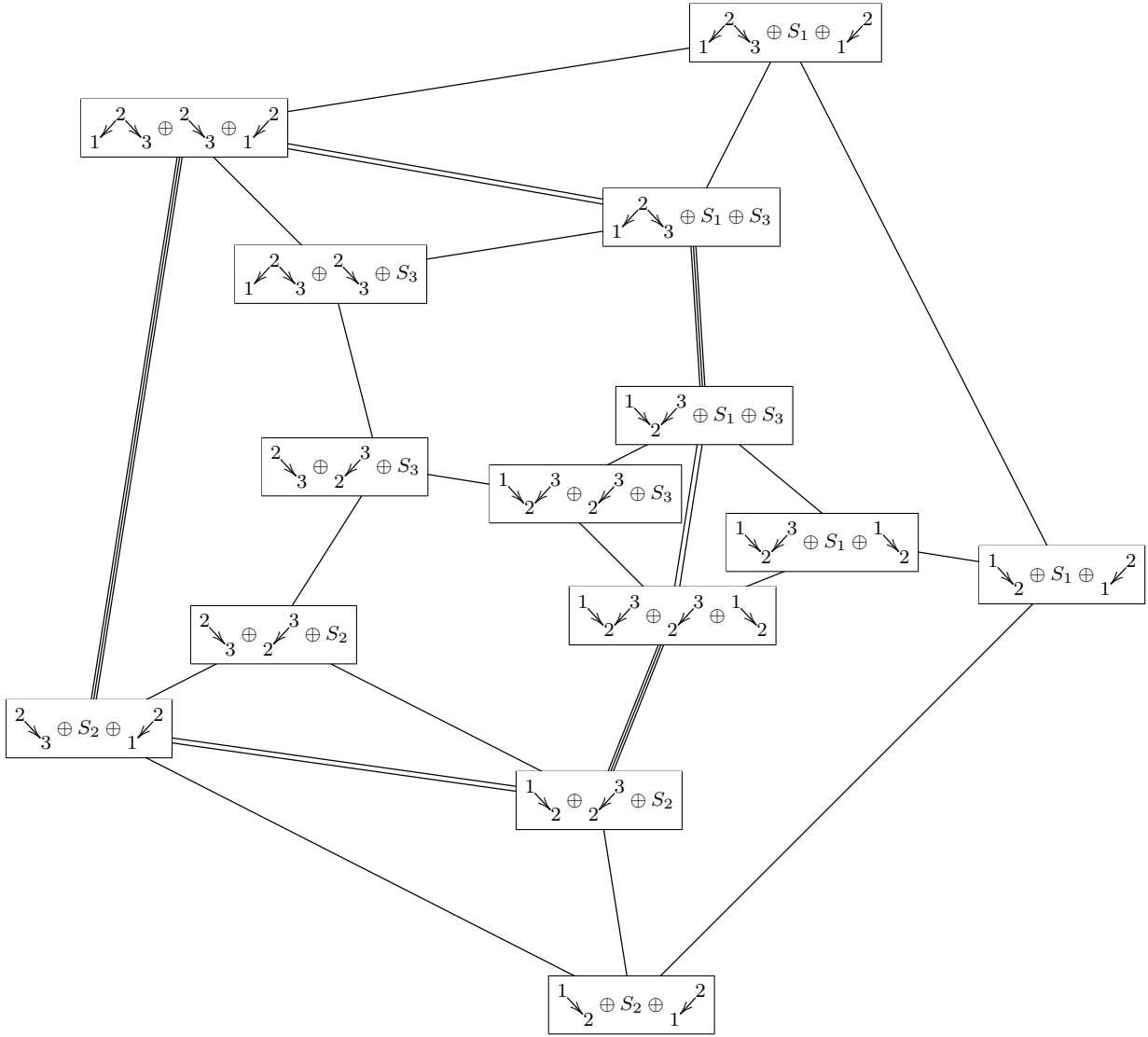
\begin{figure}
  \makeatletter
    \tiny  \objectmargin={-1pt}
    $$\begin{xy}
    		="D",0;<.58mm,.58mm>:
   		{\xypolygon5"O"{~:{(0,0):}~>{}}},
   		{\xypolygon5"A"{~:{(-30,0):}~>{}}},
    		{\xypolygon5"B"{~:{(54,0):}~>{}}},
    		{\xypolygon5"C"{~:{(90,0):}~>{}}},
    		"B1";"C1"**\dir{-},
    		"B2";"C2"**\dir{-},
    		"B3";"C3"**\dir{-},
    		"B4";"C4"**\dir{-},
    		"B5";"C5"**\dir{-},
    		"B1";"A3"**\dir{-},
    		"A3";"B5"**\dir{-},
    		"B5";"A2"**\dir{-},
    		"A2";"B4"**\dir{-},
    		"B4";"B3"**\dir{-},
    		"B3";"A5"**\dir{-},
    		"A5";"B2"**\dir{-},
    		"B2";"B1"**\dir{-},
    		"A3";"O1"**\dir{-},
    		"A2";"O1"**\dir{-},
    		"A5";"O1"**\dir{-},
    		"C1";"C2"**\dir{-},
    		"C2";"C3"**\dir{-},
    		"C3";"C4"**\dir{-},
    		"C4";"C5"**\dir{-},
    		"C5";"C1"**\dir{-},
				"A2";"B4"**\dir{=},
				"B1";"A3"**\dir{=},
 				"C2";"C3"**\dir{=},
				"A2";"A3"**\dir{=},
				"C2";"B1"**\dir{=},
				"C3";"B4"**\dir{=},
                "O1"*{\ensb{10mm}{0mm}{\utds}{\tds}{S_3}},
		"A2"*{\ensb{0mm}{0mm}{\utds}{\tds}{\uds}},
		"A3"*{\ensb{0mm}{0mm}{\utds}{S_1}{S_3}},
		"A5"*{\ensb{0mm}{10mm}{\dts}{\tds}{S_3}},
		"B1"*{\ensb{0mm}{0mm}{\duts}{S_1}{S_3}},
		"B2"*{\ensb{0mm}{0mm}{\duts}{\dts}{S_3}},
		"B3"*{\ensb{0mm}{0mm}{\dts}{\tds}{S_2}},
		"B4"*{\ensb{0mm}{0mm}{\uds}{\tds}{S_2}},
		"B5"*{\ensb{0mm}{10mm}{\utds}{S_1}{\uds}},
		"C1"*{\ensb{0mm}{0mm}{\duts}{S_1}{\dus}},
		"C2"*{\ensb{0mm}{0mm}{\duts}{\dts}{\dus}},
		"C3"*{\ensb{0mm}{0mm}{\dts}{S_2}{\dus}},
		"C4"*{\ensb{0mm}{0mm}{\uds}{S_2}{\dus}},
		"C5"*{\ensb{0mm}{0mm}{\uds}{S_1}{\dus}},
    \end{xy}$$
\makeatother
 \caption{Exchange graph of $\Z/2\Z$-stable maximal rigid objects} \label{exchnss}
\end{figure}

The exchange graph of the $\Z/2\Z$-stable basic maximal rigid objects of $\md \Pi_Q$ is presented on Figure \ref{exchnss}, in relation to the exchange graph of the basic maximal rigid objects. It permits, in view of Figure \ref{assoc} to describe the clusters of $\C[N']$:
{\small
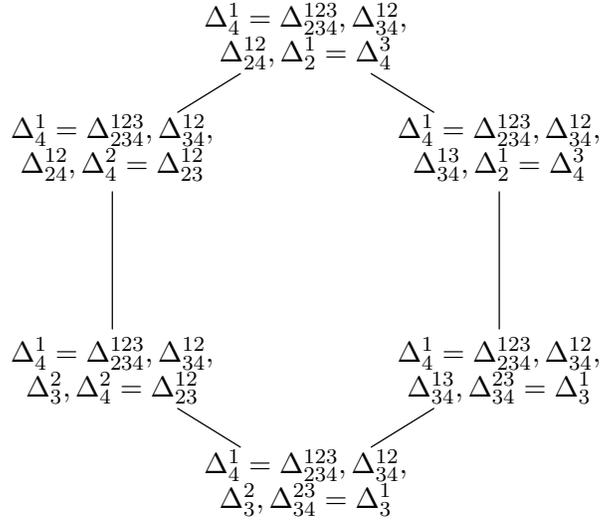
$$\begin{xy}
    		="D",0;<.7mm,.7mm>:
   		{\xypolygon6"A"{~*{\ifcase\xypolynode\or \ensc{\dta{13}{34}, \dta{23}{34} = \dta{1}{3}} \or \ensc{\dta{13}{34}, \dta{1}{2} = \dta{3}{4}} \or \ensc{\dta{12}{24}, \dta{1}{2} = \dta{3}{4}} \or \ensc{\dta{12}{24}, \dta{2}{4} = \dta{12}{23}} \or \ensc{\dta{2}{3}, \dta{2}{4} = \dta{12}{23}} \or \ensc{\dta{2}{3}, \dta{23}{34} = \dta{1}{3}} \fi }~:{(0,-30):}~={15}~>{-}}},
    \end{xy}.$$
}

\section{Scope of these results and consequences}

The example presented here can be generalized to the coordinate rings of:
\begin{itemize}
 \item The groups of the form
  $$N(w) = N \cap \left(w^{-1} N_- w \right) \quad \text{and} \quad N^w = N \cap \left(B_-w B_-\right)$$
  where $N$ is a maximal unipotent subgroup of a Kac-Moody group, $N_-$ its opposite unipotent group, $B_-$ the corresponding Borel subgroup, and $w$ is an element of the corresponding Weyl group. In particular, if $N$ is of Lie type and $w$ is the longest element, then $N(w) = N$.
 \item Partial flag varieties corresponding to classical Lie groups.
\end{itemize}

These results were obtained in \cite{GLSpart} and \cite{GLS} for the simply-laced cases and in \cite{Dem} for the non simply-laced cases.

It permits for example to prove in these cases that all the cluster monomials (products of elements of a same extended cluster) are linearly independent (result which is now generalized but was new at that time) and other more specific results (for example the classification of partial flag varieties the coordinate rings of which have finite cluster type, that is a finite number of clusters).

%%%%　参考文献の書式 %%%%%%%%%%%%%%%

%% \bibliographystyle{amsplain}

\ifx\undefined\bysame
\newcommand{\bysame}{\leavevmode\hbox to3em{\hrulefill}\,}
\fi

\end{document}